\documentclass[12pt]{amsart}
\usepackage{amsmath,amssymb}

\newtheorem{theorem}{Theorem}
\newtheorem{proposition}[theorem]{Proposition}
\newtheorem{lemma}[theorem]{Lemma}
\newtheorem{corollary}[theorem]{Corollary}
\newtheorem{definition}[theorem]{Definition}
\newtheorem{remark}[theorem]{Remark}

\def\proof{\noindent{\it Proof.}\ }
\def\QED{\hfill QED}

\def\Q{\mathcal P}
\def\Q{\mathcal Q}
\def\Z{\mathbb Z}

\begin{document}

\title[Elementary proof \& apps of Hall-Littlewood generating function]
{Elementary proof and application of the generating functions
for generalized Hall--Littlewood 
functions}

\author{Hiroshi Naruse}

\date{}


\maketitle

\begin{abstract}
In this note we define a generalization of 
Hall--Littlewood symmetric functions 
using a formal group law and we give
an elementary proof of the generating function formula
for these generalized 
functions.
We also 
present some applications
of this formula.
\end{abstract}

\section{Introduction}
Let $r$ and $n$ $(1\leq r\leq n)$ be positive integers.
Given a sequence 
$\lambda=(\lambda_1,\ldots, \lambda_r)$ of positive integers of length $r$,
a formal power series $T(x,y)$,  
and a formal Laurent series $H(z)$,
we define the 
generalized Hall--Littlewood symmetric function
$\Q^{F}_{\lambda_1,\ldots,\lambda_r}(x_1,\ldots,x_n;T,H)$
using a formal group law $F(x,y)$.
In the case of
the
additive group law 
$F(x,y)=x+y$
with $T(x,y)=x-t y$ and $H(z)=1$
this becomes the usual 
Hall--Littlewood $Q$-symmetric function 
defined in 
\cite{Mac}(cf. Definition 2 below).
Nakagawa and the author
have previously defined the universal Hall--Littlewood function $H^L_\lambda(x_1,\ldots,x_n;t)$
for a partition $\lambda$ \cite{NN}.
When $\lambda$ is a strict partition, 
definition of \cite{NN}
provides 
a special 
case of a
generalized Hall--Littlewood function,
but in general
functions are different.

We do not discuss in this paper
related
geometries
here, but we formulate the generating 
function and the proof of its properties
in almost the same way as in \cite{Mac}.
After this elementary proof, 
we found another simple proof
by the means of
the push-forward formula \cite{NN2}.
However, the proof presented in this paper only uses elementary algebraic manipulation
and does not use 
any topological arguments, such as 
Quillen's push-forward formula.
Actually,
Quillen's formula
(\cite[Theorem 1]{Qui})
can be derived from Lemma 4 below, applied to
$\ell=0, k=1$ and $T(x,y)=1$,
as well as Bressler-Evens formula
\cite[Theorem 2.5, Corollary 2.6]{NN}.

\section{
Generating function formulas and their proofs}

\subsection{Formal group law}

A formal group law $F(x,y)$ is a formal power series in two variables 
$x$ and $y$,
  $$F(x,y)=x+y+\sum_{i,j\geq1} a_{i,j}\: x^i y^j,$$
such that  $a_{i,j}$'s 
satisfy the following two properties:
 
\begin{itemize}
\item[(1)] commutativity $F(x,y)=F(y,x)$, and

\item[(2)] associativity  $F(x,F(y,z))=F(F(x,y),z)$

\end{itemize}
 (cf. \cite{LM}) .

\noindent
For a given formal group law $F$, we define an associated
 formal power series
$P_F(z)$ as follows:
$$P_F(z):=1+\sum_{i=1}^\infty a_{1,i} z^i.$$ 

\noindent
In addition, set $P_F(x,y):=\displaystyle\frac{x-y}{F(x,\bar{y})}$, where 
$\bar{y}$ is the inverse of $y$, i.e., it satisfies $F(y,\bar{y})=0$.

\begin{lemma}
$P_F(x,y)$ is invertible and $P_F(x,x)=P_F(x)$.
\end{lemma}

\proof
We use the following notation for the derivative with respect to the first variable of $F$:
$$F_1(x,y):=\partial_x F(x,y)=1+\displaystyle\sum_{i,j\geq 1} i \:a_{i,j}\: x^{i-1} y^{j}.$$
\noindent
Then, by definition, $P_F(x)=F_1(0,x)$.
We can use associativity of $F$
$$F(x,F(y,z))=F(F(x,y),z).$$
\noindent
Differentiating both sides with respect to variable $x$, we get
$$F_1(x,F(y,z))=F_1(F(x,y),z) F_1(x,y).$$
\noindent
By substituting $x=0$ and $z=\bar{y}$, we have
$$1=F_1(y,\bar{y}) F_1(0,y).$$
\noindent
Finally, using this equation, we can deduce
$$P_F(y,y)=\lim_{x\to y} \frac{x-y}{F(x,\bar{y})}=\frac{1}{F_1(y,\bar{y})}=F_1(0,y)=P_F(y),$$
\noindent
which is what we wanted to prove.

\QED

For the rest of this paper, we assume that the formal group law $F(x,y)$ is written as
$F(x,y)=\exp_F( \log_F(x)+\log_F(y))$, using appropriate 
formal power series $\exp_F$ and $\log_F$.
Then we can define formal multiplication by a variable $t$ as
$[t]_F  x:=\exp_F(t \log_F(x))$. The inverse of $x$ can be written as $\bar{x}=[-1]_F x$.
We sometimes drop the $F$ notation  and just write $[t] x$.

\subsection{Generalized Hall--Littlewood function}

We now generalize the classical Hall--Littlewood function
using the formal group law $F$.
Actually, we further generalize Hall--Littlewood function
using 
a formal power series $T(x,y)$ and
a formal Laurent series $H(z)$.
For a formal Laurent series $L(z)$, we 
use notation $[z^i] L(z)$ for the
coefficient of $z^i$ in $L(z)$.
We also allow the case of $L(z)=\sum_{i\in \Z} c_i z^{i}$
 with coefficients to be non-zero
 for all $i\in \Z$.
For example,
$\frac{1}{1+\beta z}+\frac{x}{z-x}=
\cdots+(-\beta)^2 z^2+(-\beta)z+1+(x/z)+(x/z)^2+\cdots$.
Given $T(x,y)$ and $H(z)$,
we define the generalized Hall--Littlewood function 
$\Q^{F}_{\lambda_1,\ldots,\lambda_r}(x_1,\ldots,x_n;T,H)$ 
for
a positive integer sequence $\lambda=(\lambda_1,\ldots,\lambda_r)$
as
follows. Let $S_n$ be
symmetric group acting on the variables
 $x_1,\ldots,x_n$.
\begin{definition}The generalized Hall--Littlewood function
associated with $F,T,$ and $H$ is
$$\Q^{F}_{\lambda_1,\ldots,\lambda_r}(x_1,\ldots,x_n;T,H):=
\displaystyle\sum_{w\in S_n/(S_1^r\times S_{n-r})} 
w\left( x_1^{[\lambda_1]}\cdots x_r^{[\lambda_r]}
\displaystyle  \prod_{1\leq i\leq r}\prod_{i<j\leq n} \frac{T(x_i, {x}_j)}{F(x_i, \bar{x}_j)}\right)
,$$
\noindent
where $x^{[k]}:=H(x)T(x, {x}) x^{k-1}$ for any positive integer $k$.
\end{definition}
In the following, we sometimes drop the $F,T,H$ notations 
and simply write this function as
$\Q_\lambda(x_1,\ldots, x_n)$.

\begin{definition}
We focus on the following two cases:

\begin{enumerate}

\item

$T(x,y)=F(x,[t]\bar{y})$, $H(z)=1$ and define

$HQ_\lambda(x;t):=$

\hspace{1cm}$
\displaystyle\sum_{w\in S_n/(S_1^r\times S_{n-r})} 
w\left( F(x_1,[t]\bar{x}_1)x_1^{\lambda_1-1}\cdots F(x_r,[t]\bar{x}_r)x_r^{\lambda_r-1}
\displaystyle  \prod_{1\leq i\leq r}\prod_{i<j\leq n} \frac{F(x_i, [t]\bar{x}_j)}{F(x_i, \bar{x}_j)}\right);$
\\

\item

$T(x,y)=F(x,[t]\bar{y})$, $H(z)=\frac{z}{T(z,z)}$
and define

$HP_\lambda(x;t):=
\displaystyle\sum_{w\in S_n/(S_1^r\times S_{n-r})} 
w\left( x_1^{\lambda_1}\cdots x_r^{\lambda_r}
\displaystyle  \prod_{1\leq i\leq r}\prod_{i<j\leq n} \frac{F(x_i, [t]\bar{x}_j)}{F(x_i, \bar{x}_j)}\right)
.$

\end{enumerate}

\end{definition}
When $F(x,y)=x+y$, 
$HQ_\lambda(x;t)$ becomes the Hall--Littlewood $Q$-function 
and for a strict partiton $\lambda$, 
$HP_\lambda(x;t)$ becomes  the Hall--Littlewood $P$-function 
defined in \cite{Mac}.
For the connective $K$-theory case
$F=F_{CK}(x,y)=x+y+\beta x y$
with $t=-1$,
they are the $K$-theoretic Schur
$Q$- and $P$-functions $GQ_\lambda(x)$ and $GP_\lambda(x)$ 
respectively
defined in \cite{IN}.

\begin{lemma}
Let
$T(x,y)$ be a formal power series in $x$ and $y$,
and let
$H(z)$ be a formal Laurent series in $z$ with
a pole at $z=0$ of order  $\ell\geq 0$,
i.e.,
$H(z)=\sum_{i\geq -\ell} c_i z^i$ with $c_{-\ell}\neq 0$.
Then, for $k> \ell$, we have
 $$[z^{-k}] 
 \left(\frac{H(z)}{P_F(z)}\prod_{i=1}^n \frac{T(z, {x}_i)}{F(z, \bar{x}_i)}
 \right)
 =
\sum_{i=1}^{n} {x_i^{k-1} H(x_i)T(x_i, {x}_i)}
\prod_{j\neq i}\frac{T(x_i, {x}_j)}{F(x_i, \bar{x}_j)}.
$$
 
\end{lemma}

\proof
The expression under consideration can be written as
$$\frac{H(z)}{P_F(z)}\prod_{i=1}^n \frac{T(z, {x}_i)}{F(z, \bar{x}_i)}
=
\frac{H(z)}{P_F(z)} \prod_{i=1}^n \frac{P_F(z,{x}_i) T(z, {x}_i) }{z-{x}_i}
.$$

\noindent
We define
$$R_F(z):=
z^{\ell}
\frac{H(z)}{P_F(z)}
\prod_{i=1}^{n} \left(P_F(z,x_i) T(z, {x}_i)\right)
-
\sum_{i=1}^{n}
x_i^\ell H(x_i) T(x_i, {x}_i) \prod_{j\neq i}\frac{T(x_i, {x}_j)(z-x_j)}{F(x_i, \bar{x}_j)}.
$$
\noindent
This is a formal power series in variable $z$.
When specializing $z=x_k$  ($k=1,2,\ldots, n$), it vanishes.
Therefore, we can write
$$R_F(z)=B(z) \displaystyle \prod_{i=1}^{n} (z-x_i) ,$$
where $B(z)$ is a formal power series in $z$.
By dividing this by $\displaystyle\prod_{i=1}^{n} (z-x_i)$, we get

$$z^{\ell}
\frac{H(z)}{P_F(z)}
\prod_{i=1}^{n} \frac{P_F(z,x_i) T(z, {x}_i)}{z-x_i}
=
B(z)+
\sum_{i=1}^{n}
\frac{
x_i^\ell H(x_i) T(x_i, {x}_i)}{z-x_i} \prod_{j\neq i}\frac{T(x_i, {x}_j)}{F(x_i, \bar{x}_j)}.
$$
\noindent
By taking the coefficients of $z^{\ell-k}$ in both sides,
we
obtain 
the claim of the lemma.
\QED

\vspace{0.5cm}

\begin{corollary}
Setting
$T(x,y)=F(x,[t]\bar{y})$
and
$H(z)=1$,
we get the equality
$$
\frac{1}{P_F(z)}\prod_{i=1}^n \frac{F(z, [t] \bar{x}_i)}{F(z, \bar{x}_i)}
=
N_F(z;x_1,\ldots,x_n;t)
+
\sum_{i=1}^{n} \frac{F(x_i, [t] \bar{x}_i)}{z-x_i} 
\prod_{j\neq i}\frac{F(x_i, [t] \bar{x}_j)}{F(x_i, \bar{x}_j)},
$$
\noindent
where $N_F(z;x_1,\ldots,x_n;t)$ is a formal power series in 
variable $z$.

\end{corollary}

\begin{remark}
In the connective $K$-theory case $F_{CK}(x,y)=x+y+\beta x y$
 we can prove that
$$N_{F_{CK}}(z;x_1,\ldots,x_n;t)=
\frac{1}{1+\beta z},$$
(cf. Proposition 12).

\noindent
In the elliptic cohomology case $F_{Ell}(x,y)=\frac{x+y+\beta x y}{1-\gamma x y}$ we have

$$N_{F_{Ell}}(z;x_1;t)=\frac{1}{1+\beta z+\gamma z^2}+
\frac{\gamma w_1(x_1+w_1+\beta x_1 w_1)}{(1-\gamma  w_1 x_1)(1-\gamma w_1 z)},$$
where $w_1=[t] \bar{x}_1$.

\end{remark}

\begin{corollary}
Setting
$T(x,y)=F(x,[t]\bar{y})$ and $H(z)=\frac{z}{T(z,z)}$,
we
obtain the equality
$$
\frac{z}{F(z, [t] \bar{z})}\frac{1}{P_F(z)}\prod_{i=1}^n \frac{F(z, [t] \bar{x}_i)}{F(z, \bar{x}_i)}
=
M_F(z;x_1,\ldots,x_n;t)
+
\sum_{i=1}^{n} \frac{x_i}{z-x_i} 
\prod_{j\neq i}\frac{F(x_i, [t] \bar{x}_j)}{F(x_i, \bar{x}_j)},
$$
\noindent
where $M_F(z;x_1,\ldots,x_n;t)$ is a formal power series in 
variable $z$.

\end{corollary}

\begin{remark}
In the connective K-theory case $F_{CK}(x,y)=x+y+\beta x y$
 we can prove that
$$M_{F_{CK}}(z;x_1,\ldots,x_n;t=-1)=\frac{1}{1+\beta z}
-\frac{\prod_{i=1}^{n}(1+\beta x_i)}{2+\beta z}
$$
(cf. Proposition 13).
\end{remark}

\begin{theorem}

Assume that
 $T(x,y)$ is a formal power series
and $T(x,0)=x$.
Assume further that 
$H(z)$ is a formal Laurent series in $z$ with a pole at zero
 of order  $\ell\geq 0$.

Let $A(z):=\displaystyle\frac{H(z)}{P_F(z)}
\prod_{i=1}^{n}\frac{T(z, {x}_i)}{F(z, \bar{x}_i)}$.
Set
$$A(z_1,\ldots,z_r):=\prod_{i=1}^{r} A(z_i) 
\prod_{1\leq i<j\leq r}\frac{F(z_j, \bar{z}_i)}{T(z_j, {z}_i)}.$$
\noindent
Then
for a sequence of positive integers $\lambda_1,\ldots,\lambda_r>\ell$ , $(r\leq n)$, we have 
 $$[z_1^{-\lambda_1}\cdots z_r^{-\lambda_r}]A(z_1,\ldots,z_r)=\Q^F_\lambda(x_1,\ldots, x_n;T,H).$$
\end{theorem}

\proof
(See
\cite{Mac} p.211  proof of (2.15)).
We will show that for a sequence of positive integers $\lambda=(\lambda_1,\ldots,\lambda_r)$
the coefficient of 
$z^{-\lambda}:=z_1^{-\lambda_1}\cdots z_r^{-\lambda_r}$ 
in $A(z_1,\ldots,z_r)$ is
$\Q_\lambda(x_1,\ldots, x_n)$.
For the case $r=1$ it follows by Lemma 4.
\noindent
Assume $r\geq 2$.
By definition, 
$$A(z_1,\ldots,z_r)=A(z_2,\ldots,z_r)\times A(z_1) 
\prod_{j=2}^r\frac{F(z_j,\bar{z}_1)}{T(z_j, {z}_1)}.$$
We can expand  the
product
in the formula above
 as a formal power series in $z_1$ 
as follows.
$$\prod_{j=2}^r\frac{F(z_j,\bar{z}_1)}{T(z_j, {z}_1)}
=\sum_{m\geq 0} f_m(z_2,\ldots,z_r) z_1^m.
$$
On the other hand, by Lemma 4
the coefficient of $z_1^{-\lambda_1-m}$ in
$A(z_1)$ is $\Q_{\lambda_1+m}(x_1,\ldots,x_n)$.
Therefore, the cofefficient of $z_1^{-\lambda_1}$ in $A(z_1,\ldots,z_r)$
is 
$$[z_1^{-\lambda_1}]A(z_1,\ldots,z_r)=
A(z_2,\ldots,z_r)\times\left(
\sum_{m\geq 0} f_m(z_2,\ldots,z_r) \Q_{\lambda_1+m}(x_1,\ldots,x_n) 
\right).
$$
However, by definition,
$$\Q_{\lambda_1+m}(x_1,\ldots,x_n)=
\sum_{i=1}^{n} x_i^{\lambda_1+m}H(x_i)\frac{T(x_i, {x}_i)}{x_i}
\prod_{1\leq j\leq n, j\neq i}\frac{T(x_i, {x}_j)}{F(x_i,\bar{x}_j)}.
$$
\noindent
Therefore, we can substitute this into the previous equation 
to obtain

$\begin{array}{ccl}
&&[z_1^{-\lambda_1}]A(z_1,\ldots,z_r)\\
&=&
A(z_2,\ldots,z_r)\times\left(\displaystyle
\sum_{m\geq 0}
f_m(z_2,\ldots,z_r)
\sum_{i=1}^{n} x_i^{\lambda_1+m}H(x_i)\frac{T(x_i, {x}_i)}{x_i}
\prod_{1\leq j\leq n, j\neq i}\frac{T(x_i, {x}_j)}{F(x_i,\bar{x}_j)} 
\right)\\[0.8cm]
&=&
A(z_2,\ldots,z_r)\times\left(\displaystyle
\sum_{i=1}^{n} 
\left( \sum_{m\geq 0} x_i^{m}f_m(z_2,\ldots,z_r) \right)
x_i^{\lambda_1} H(x_i)\frac{T(x_i,{x}_i)}{x_i}
\prod_{1\leq j\leq n, j\neq i}\frac{T(x_i,{x}_j)}{F(x_i,\bar{x}_j)} 
\right)\\[0.8cm]
&=&\displaystyle
\sum_{i=1}^{n} 
\left(A(z_2,\ldots,z_r) \prod_{j=2}^r\frac{F(z_j,\bar{x}_i)}{T(z_j, {x}_i)} \right)
x_i^{\lambda_1} H(x_i)\frac{T(x_i, {x}_i)}{x_i}
\prod_{1\leq j\leq n, j\neq i}\frac{T(x_i, {x}_j)}{F(x_i,\bar{x}_j)}\\[0.8cm]
&=&\displaystyle
\sum_{i=1}^{n} 
\left(A^{(i)}(z_2,\ldots,z_r) \right)
x_i^{\lambda_1} H(x_i)\frac{T(x_i, {x}_i)}{x_i}
\prod_{1\leq j\leq n, j\neq i}\frac{T(x_i, {x}_j)}{F(x_i,\bar{x}_j)} ,\\
\end{array}
$
where we denoted by $A^{(i)}(z_2,\ldots,z_r)$  the result of setting $x_i=0$
in $A(z_2,\ldots,z_r)$.
Here we used the property $T(x,0)=x$.
\noindent
By the induction hypothesis, 
the coefficient of $z_2^{-\lambda_2}\cdots z_r^{-\lambda_r}$ in 
$A(z_2,\ldots,z_r)$
is 
$\Q_{(\lambda_2,\ldots,\lambda_r)}(x_1,\ldots,x_n)$;
therefore, the coefficient of $z_1^{-\lambda_1}z_2^{-\lambda_2}\cdots z_r^{-\lambda_r}$ 
in $A(z_1,\ldots,z_r)$ is
$$
\sum_{i=1}^{n} 
\left(\Q^{(i)}_{(\lambda_2,\ldots,\lambda_r)}(x_1,\ldots,x_n) \right)
x_i^{\lambda_1} H(x_i)\frac{T(x_i, {x}_i)}{x_i}
\prod_{1\leq j\leq n, j\neq i}\frac{T(x_i, {x}_j)}{F(x_i,\bar{x}_j)}, 
$$
which can  easily be seen to equal 
$\Q_{(\lambda_1,\ldots,\lambda_r)}(x_1,\ldots,x_n).$
\QED

\begin{remark}

We can modify the above proof so that
each $A(z)$ has a different $H(z)$. This gives us a formula 
for the factorial
version and its determinant-Pfaffian formula.
For example, for a given 
sequence
$\lambda=(\lambda_1,\ldots,\lambda_r)$
of positive integers, we have

$
[z^{-\lambda}]
\left(
A_1^{(\lambda_1-1)}(z_1)\cdots A_r^{(\lambda_r-1)}(z_r)
\displaystyle\prod_{1\leq i<j\leq r}\frac{F(z_j, \bar{z}_i)}{T(z_j, {z}_i)}
\right)
$

$$
=\displaystyle\sum_{w\in S_n/(S_1^r\times S_{n-r})} 
w\left( (x_1|b)^{[\lambda_1]}\cdots (x_r|b)^{[\lambda_r]}
\displaystyle  \prod_{1\leq i\leq r}\prod_{i<j\leq n} \frac{T(x_i, {x}_j)}{F(x_i, \bar{x}_j)}\right),
$$
where
$A_i^{(k)}(z):=\displaystyle\frac{H_i^{(k)}(z)}{P_F(z)}
\prod_{i=1}^{n}\frac{T(z, {x}_i)}{F(z, \bar{x}_i)}$
,
if we set
$H_i^{(k)}(z)=(z|b)^{k}/z^{k}$, $(z|b)^{k}=\displaystyle\prod_{i=1}^{k} F(z,b_i)$
as the factorial power
and  $(x|b)^{[k]}:=T(x, {x}) (x|b)^{k-1}$
(cf. similar formulas in
\cite{NN2}).

\end{remark}

We set
$x\oplus y:=F_{CK}(x,{y})=x+y+\beta x y$ and
$x\ominus y:=F_{CK}(x,\bar{y})=\frac{x-y}{1+\beta y}$.
A geometric proof of Corollary 11 below is given in \cite{HIMN}.

\begin{corollary}(Determinant-Pfaffian formula)
Consider the connective $K$-theory case $F=F_{CK}(x,y)$
and assume $\lambda$
is a partition of length $r$,
$T(x,y)=F(x,[t]\bar{y})$,
and
$H(z)=1$ or $H(z)=\frac{z}{T(z,z)}$.
Then
for $t=0$,
 we have the determinantal formulas
$$\begin{array}{lll}
\Q_{\lambda_1,\ldots,\lambda_r}(x_1,\ldots,x_n)
&=&\det\left( [z_i^{-(\lambda_i+j-i)}] \displaystyle\frac{1}{(1+\beta z_i)^{r-i}}
A(z_i)\right)_{r\times r}\\[0.3cm]
&=&\det\left( [z_i^{-(\lambda_i+j-i)}] \displaystyle\frac{1}{(1+\beta z_i)^{j-i}}
A(z_i)\right)_{r\times r},
\end{array}
$$
and for $t=-1$, if we assume $r$ is even, we have the
Pfaffian formula
$$
\Q_{\lambda_1,\ldots,\lambda_r}(x_1,\ldots,x_n)={\rm Pf} \left( 
[z_i^{-\lambda_i} z_j^{-\lambda_j}]
\frac{1}{(1+\beta z_i)^{r-i-1}}\frac{1}{(1+\beta z_i)^{r-j}}
  A(z_i)A(z_j) \frac{z_j\ominus z_i}{z_j\oplus z_i}\right)_{r\times r}.
$$

\end{corollary}

\proof
For the determinantal formulas, 
using the identity

$$
\frac{z_j\ominus z_i}{z_j}=\frac{1-z_i/z_j}{1+\beta z_i}=1-\frac{\bar{z}_i}{\bar{z}_j},
$$
the assertion follows from Theorem 9.
For the Pfaffian formula,  by Theorem 9,
we have

$\Q_\lambda(x)
=
[z^{-\lambda}]
A(z_1)\cdots A(z_{r})
\displaystyle
\prod_{1\leq i<j\leq r}
\frac{z_j\ominus z_i}{z_j\oplus z_i}\\
=
[z^{-\lambda}]
A(z_1)\cdots A(z_{r})
\displaystyle
\prod_{i=1}^{r-1} \frac{1}{(1+\beta z_i)^{r-i}}
\prod_{1\leq i<j\leq r}
\frac{z_j- z_i}{z_j\oplus z_i}\\
=
[z^{-\lambda}]
A(z_1)\cdots A(z_{r})
\displaystyle
\prod_{i=1}^{r-1} \frac{1}{(1+\beta z_i)^{r-i}}
{\rm Pf} \left(\frac{z_j- z_i}{z_j\oplus z_i}\right)
\\
=[z^{-\lambda}]
\displaystyle
{\rm Pf} \left(
A(z_i) A(z_{j})
\frac{1}{(1+\beta z_i)^{r-i}}
\frac{1}{(1+\beta z_j)^{r-j}}
\frac{z_j- z_i}{z_j\oplus z_i}\right)
\\
=
\displaystyle
{\rm Pf} \left([z_i^{-\lambda_i}z_j^{-\lambda_j}]
A(z_i) A(z_{j})
\frac{1}{(1+\beta z_i)^{r-i-1}}
\frac{1}{(1+\beta z_j)^{r-j}}
\frac{z_j\ominus z_i}{z_j\oplus z_i}\right).
$

\QED

\subsection{Non-negative power part}

\begin{proposition}
For $F=F_{CK}(x,y)=x+y+\beta x y,$ we have

$$N_F(z;x_1,\ldots,x_n;t)=\frac{1}{1+\beta z}.$$

\end{proposition}

\proof
We will prove
$$
\frac{1}{1+\beta z}\prod_{i=1}^n \frac{z\ominus w_i}{z\ominus{x}_i}
=
\frac{1}{1+\beta z}+
\sum_{i=1}^{n} \frac{x_i\ominus w_i}{z-x_i} 
\prod_{j\neq i}\frac{x_i\ominus w_j}{x_i\ominus {x}_j}.
$$
Then we can substitute $w_i=[t] x_i$ ($i=1,\ldots,n$)  to get the assertion of the proposition.
Multiplying
by
$(1+\beta z)\prod_{i=1}^{n}(z\ominus x_i)$, we must prove 
$$\prod_{i=1}^n {(z\ominus w_i)}
=
\prod_{i=1}^n {(z\ominus x_i)}
+
(1+\beta z)
\sum_{i=1}^{n} 
\frac{1}{1+\beta x_i}
\prod_{j=1}^n {(x_i\ominus w_j)} 
\prod_{j\neq i}\frac{z\ominus {x}_j}{x_i\ominus {x}_j}.
$$

\noindent
But both sides are polynomial in $z$ of degree $n$, and the equation has  $n$ solutions
$z=x_i$  ($i=1,\ldots,n$).
So we need to prove that the top degree terms have the same coefficient.
It gives the equality
$$
\frac{\prod_{i=1}^{n} (1+\beta x_i)}{\prod_{i=1}^{n} (1+\beta w_i)}=
1
+\beta
\sum_{i=1}^{n}
\frac{\prod_{j=1}^{n} (x_i\ominus w_j)}{\prod_{j\neq i} (x_i\ominus x_j)}
.$$

\QED

\begin{proposition}

$t=-1$ and $F=F_{CK}$ is the connective $K$-theory case, and 
we have

$$M_F(z;x_1,\ldots,x_n;t=-1 )=
\frac{1}{(1+\beta z)(2+\beta z)}
-\frac{\beta GP_1(x_1,\ldots,x_n)}{2+\beta z}.$$

Note that 
$ 1+\beta GP_1(x_1,\ldots,x_n)=\displaystyle\prod_{i=1}^{n}(1+\beta x_i)$.
\end{proposition}

\proof
We need to show
$$
\frac{1}{1+\beta z}\frac{1}{2+\beta z}\prod_{i=1}^{n}\frac{z\oplus x_i}{z\ominus x_i}
=
\frac{1}{1+\beta z}-\frac{1+\beta GP_1}{2+\beta z}
+
\sum_{i=1}^{n} \frac{x_i}{z-x_i}\prod_{j\neq i}\frac{x_i\oplus x_j}{x_i\ominus x_j}.
$$

\noindent
For this we multiply by 
$(1+\beta z)(2+\beta z)\displaystyle\prod_{i=1}^{n} (z\ominus x_i)$
both sides.

\noindent
Therefore, we need to show

$\displaystyle
\prod_{i=1}^{n}{(z\oplus x_i)}
=
\left((2+\beta z)-(1+\beta z)\prod_{i=1}^{n}(1+\beta x_i)\right)
\prod_{i=1}^{n}{(z\ominus x_i)}$

$\hspace{3cm}\displaystyle
+(1+\beta z)(2+\beta z)\sum_{i=1}^{n} 
\frac{z_i}{1+\beta x_i}
\prod_{j\neq i}\left(\frac{x_i\oplus x_j}{x_i\ominus x_j}(z\ominus x_j)\right).
$

\noindent
This holds for $z=x_1,\ldots,x_n$. So we need to check that the
coefficient of $z^{n+1}$ on the right hand side  is zero, 
and that the coefficients of $z^n$ are the same for both sides.
The coefficient of $z^{n+1}$ on the right hand side is
$$\beta^2 \sum_{i=1}^{n} \frac{x_i}{1+\beta x_i}\prod_{j\neq i} \frac{x_i\oplus x_j}{x_i\ominus x_j}
\frac{1}{1+\beta x_j}-\beta^2GP_1=0.
$$

\noindent
The coefficient of $z^{n}$ on the right hand side is

\noindent
$\displaystyle\frac{3\beta {GP_1}+\beta^2(GP_2-(x_1+\cdots+x_n)GP_1)}{1+\beta GP_1}+
\frac{2-(1+\beta GP_1)+(x_1+\cdots+x_n)\beta^2 GP_1}{1+\beta GP_1}\\
=1+\beta GP_1
.$

\noindent
As the left hand side is
$\displaystyle
\prod_{i=1}^n (z(1+\beta x_i)+x_i)
$,
the coefficient of $z^n$ is
$$
\prod_{i=1}^n (1+\beta x_i)
=1+\beta GP_1,
$$
so we get the result.

\QED

\section{Application}

In this section we explain
how the generating function can be used to
give a formula such as the Pieri rule.

\subsection{Two row in terms of one rows}

We explain 
the case of $K$-theory and $t=-1$, i.e.
we express $GQ_{k,\ell}$ in terms of $GQ_p\times GQ_q$.

\begin{lemma}We have the identity
$$\frac{z_2\ominus z_1}{z_2\oplus z_1}
=\sum_{i\geq j\geq 0} g_{i,j} z_1^i z_2^{-j},$$
\noindent
where $g_{i,j}:=(-1)^i \beta^{i-j} \left(2\binom{i}{j}+\binom{i}{j+1}-\delta_{j,0}  \right)$
and $\delta_{i,j}$ is the Kronecker delta.

\end{lemma}

\proof
The assertion follows from the following equality:
$$\frac{z_2\ominus z_1}{z_2\oplus z_1}=\frac{(1-z_1/z_2)}{(1+\beta z_1)(1+(1/z_2+\beta)z_1)}
=\frac{2+\beta z_2}{1+(1/z_2+\beta)z_1}-\frac{1+\beta z_2}{1+\beta z_1}.$$
\QED

\begin{proposition} For $k>\ell>0$
 we have
$$GQ_{k,\ell}=\sum_{j=0}^{\ell}
\sum_{i=j}^\infty c_{i,j} GQ_{k+i} GQ_{\ell-j},$$
where 
$$
c_{i,j}=\begin{cases}
g_{i,j}& (0\leq j<\ell)\\
(-1)^{i}\beta^{i-j}\left( 2 \binom{i-1}{j-1}+\binom{i-1}{j}\right)& (j=\ell).
\end{cases}
$$

\end{proposition}

\proof
We will compare the coefficients of $z_1^{-k} z_2^{-\ell}$ on both side of 
$$A(z_1,z_2)=A(z_1) A(z_2)\frac{z_2\ominus z_1}{z_2\oplus z_1}.$$

\noindent
We expand
$A(z)=\sum_{m\in \Z} {\bf GQ}_{m} z^{m}$. Then
if $m>0$, ${\bf GQ}_{m} =GQ_m$, and if
$m\leq 0$,
 then ${\bf GQ}_{m} =(-\beta)^{-m}$ by Proposition 12.
By Lemma 14, 
it follows that
$$
c_{i,j}=\begin{cases}
g_{i,j}& (0\leq j<\ell)\\[0.3cm]
\displaystyle\sum_{i\geq s\geq j} g_{i,s} (-\beta)^{s-j}& (j=\ell).
\end{cases}
$$

Then the following binomial 
equality 
enables us to achieve the desiered result:
$$
\sum_{i\geq s\geq j} \binom{i}{s} (-1)^{s-j}=
\binom{i-1}{j-1}
\text{ for }
i\geq j>0.
$$
Note that
this equation
is derived
from the comparison of coefficients of $x^j$ in
$$
(1+x)^i \frac{1}{1+x^{-1}}=(1+x)^{i-1} x.$$

\begin{proposition} For $k>\ell>0$, we have
$$GP_{k,\ell}=\sum_{j=0}^{\ell}\sum_{i=j}^\infty d_{i,j} GP_{k+i} GP_{\ell-j},$$
where if $\ell\geq 2$
$$
d_{i,j}=\begin{cases}
g_{i,j}& (0\leq j<\ell-1)\\
(-1)^i \beta^{i-j} 2\binom{i}{j}& (j=\ell-1)\\
(-1)^{i}\beta^{i-j}\binom{i-1}{j-1}& (j=\ell),
\end{cases}
$$
and if $\ell=1$
$$
d_{i,j}=(-1)^{i}\beta^{i-j}.
$$

\end{proposition}

\proof
By Proposition 13, we have
$$
M_F(z;x_1,\ldots,x_n;t=-1 )
=\sum_{s=0}^\infty (-\beta)^s\left(
\left(1-\frac{1}{2^{s+1}}\right)-\frac{\beta GP_1(x_1,\ldots,x_n)}{2^{s+1}} \right) z^s.
$$
So, we must consider the case where $j=\ell-1,\ell$.
If $\ell\geq 2$,

$d_{i,\ell-1}=\displaystyle
g_{i,\ell-1}+\sum_{s=0}^{i-\ell} g_{i,\ell+s} 
\left(\frac{-\beta}{2}\right)^{s+1}
=(-1)^i \beta^{i-(\ell-1)}2\binom{i}{\ell-1},
$

and

$d_{i,\ell}=\displaystyle
\sum_{s=0}^{i-\ell} g_{i,\ell+s} 
\left(1-\frac{1}{2^{s+1}}\right)(-\beta)^{s}
=(-1)^i \beta^{i-\ell}\binom{i-1}{\ell-1}.
$

\noindent
If $\ell=1$,
we 
find that

$\begin{array}{ccl}
d_{i,0}&=&g_{i,0}+\displaystyle
\sum_{s=0}^{i-1} \left(\frac{-\beta}{2}\right)^{s+1} g_{i,s+1}\\
&=&(-\beta)^i \left((i+1)
-\frac{1}{2}(2\binom{i}{1}+\binom{i}{2})
+\frac{1}{2^2}(2\binom{i}{2}+\binom{i}{3})-\cdots
\right)\\
&=&
(-\beta)^{i},\\
\text{ and}\\
d_{i,1}&=&
\displaystyle\sum_{s=0}^{i-1} (-\beta)^s 
\left(1-\frac{1}{2^{s+1}}\right) g_{i,s+1}\\
&=&
(-1)^i \beta^{i-1} 
\left( (1-\frac{1}{2})(2\binom{i}{1}+\binom{i}{2})
-
(1-\frac{1}{2^2})(2\binom{i}{2}+\binom{i}{3})+
\cdots
\right)\\
&=&(-1)^i\beta^{i-1} \left(2\binom{i-1}{0}+\binom{i-1}{1}
-\frac{1}{2}(2\binom{i}{1}+\binom{i}{2})
+\frac{1}{2^2}(2\binom{i}{2}+\binom{i}{3})-\cdots
\right)\\
&=&(-1)^i\beta^{i-1}.\\
\end{array}
$

\QED

\subsection{Pieri rule for
multiplying one row by another}
For $F=F_{CK}(x,y)=x+y+\beta x y$,
using that
$$
A(z_1)A(z_2)=A(z_1,z_2)\frac{z_2\oplus z_1}{z_2\ominus z_1}
$$
 and Lemma 17, we can obtain the following Pieri rule below (Proposition 18).

\begin{lemma} We have the following equality:
$$
\frac{z_2\oplus z_1}{z_2\ominus z_1}
=(1+\beta z_1)^2+(1+\beta z_1)(2+\beta z_1)((z_1/z_2)+
(z_1/z_2)^2+(z_1/z_2)^3+\cdots).
$$
\end{lemma}
\proof
The assertion can be shown 
in a straightforward way
as follows.
$$
\frac{z_2\oplus z_1}{z_2\ominus z_1}=
\frac{(z_2+z_1+\beta z_2 z_1)(1+\beta z_1)}{z_2-z_1}
=
\frac{((1+\beta z_1)+(z_1/z_2))(1+\beta z_1)}{1-z_1/z_2}
$$

$
=(1+\beta z_1)^2+(2+\beta z_1)(1+\beta z_1)((z_1/z_2)+
(z_1/z_2)^2+(z_1/z_2)^3+\cdots)
$

\QED

\begin{proposition}
 For $k\geq l>0$, we have the following:\\

$\begin{array}{lll}
(1)\;GQ_k \times GQ_l&=&(GQ_{k,l}+2\beta GQ_{k+1,l}+\beta^2 GQ_{k+2,l})\\
&&+\displaystyle\sum_{i=1}^{l-1}
(2GQ_{k+i,l-i}+3\beta GQ_{k+1+i,l-i}+\beta^2 GQ_{k+2+i,l-i})\\
&&+(2GQ_{k+l,0}+\beta GQ_{k+l+1,0}),
\end{array}
$

$(2)\;GP_k \times GP_1=
GP_{k,1}+\beta GP_{k+1,1}+GP_{k+1}.
$

For $k\geq l\geq 2$,

\;\;$\begin{array}{lll}
GP_k \times GP_l&=&(GP_{k,l}+2\beta GP_{k+1,l}+\beta^2 GP_{k+2,l})\\
&&+\displaystyle\sum_{i=1}^{l-2}
(2GP_{k+i,l-i}+3\beta GP_{k+1+i,l-i}+\beta^2 GP_{k+2+i,l-i})\\
&&+(2GP_{k+l,1}+2\beta GP_{k+l+1,1})\\
&&+GP_{k+l}.
\end{array}
$

\end{proposition}

\proof
These follow directly from Lemma 17.
\QED

\vspace{0.5cm}

\noindent
{\bf Acknowledgements}

\vspace{0.5cm}

The author is gratefull to 
Natasha Rozhkovskaya for checking and correcting
English expressions to make this document readable.
This work was partially supported by
the Grant-in-Aid for Scientific Research (B) 16H03921,
Japan Society for the Promotion of Science.

\vspace{1cm}

Graduate School of Education, University of Yamanashi,
4-4-37, Takeda, Kofu,\\
Yamanashi, 400-8510, JAPAN

hnaruse@yamanashi.ac.jp

\end{document}